\documentclass[12pt]{article}
\usepackage{amsmath,amssymb,amsthm,amscd}
\usepackage{algorithm,algorithmic}
\usepackage[colorlinks, linkcolor=blue, anchorcolor=gree, citecolor=red]{hyperref}
\usepackage[numbers,sort&compress]{natbib}
\usepackage{graphicx}

\begin{document}

\newcommand{\Cyc}{{\rm{Cyc}}}\newcommand{\diam}{{\rm{diam}}}

\newtheorem{thm}{Theorem}[section]
\newtheorem{pro}[thm]{Proposition}
\newtheorem{lem}[thm]{Lemma}
\newtheorem{fac}[thm]{Fact}
\newtheorem{ob}[thm]{Observtion}
\newtheorem{cor}[thm]{Corollary}
\theoremstyle{definition}
\newtheorem{ex}[thm]{Example}
\newtheorem{df}[thm]{Definition}
\newtheorem{remark}[thm]{Remark}
\newcommand{\bth}{\begin{thm}}
\renewcommand{\eth}{\end{thm}}
\newcommand{\bex}{\begin{examp}}
\newcommand{\eex}{\end{examp}}
\newcommand{\bre}{\begin{remark}}
\newcommand{\ere}{\end{remark}}

\newcommand{\bal}{\begin{aligned}}
\newcommand{\eal}{\end{aligned}}
\newcommand{\beq}{\begin{equation}}
\newcommand{\eeq}{\end{equation}}
\newcommand{\ben}{\begin{equation*}}
\newcommand{\een}{\end{equation*}}

\newcommand{\bpf}{\begin{proof}}
\newcommand{\epf}{\end{proof}}
\renewcommand{\thefootnote}{}

\def\beql#1{\begin{equation}\label{#1}}
\title{\Large\bf The power index of a graph}

\author{{\sc Xuanlong Ma$^{1}$~~Min Feng$^{2,}$\footnote{Corresponding author.}~~Kaishun Wang$^1$
}\\[15pt]
{\small\em $^1$Sch. Math. Sci. {\rm \&} Lab. Math. Com. Sys., Beijing Normal University, Beijing, 100875, China}\\
{\small\em $^2$School of Science, Nanjing University of Science and Technology, Nanjing, 210094, China}\\
}

 \date{}

\maketitle

\begin{abstract}
The {\em power index} $\Theta(\Gamma)$ of a graph $\Gamma$
is the least  order of a group $G$ such that  $\Gamma$ can embed into the power graph of $G$. Furthermore, this group $G$ is  {\em $\Gamma$-optimal} if $G$ has order $\Theta(\Gamma)$.
We say that $\Gamma$ is {\em power-critical} if its order equals to $\Theta(\Gamma)$.
This paper focuses on the power indices of complete graphs, complete bipartite graphs and $1$-factors. We classify all power-critical graphs $\Gamma'$ in these three families, and  give a necessary and sufficient condition for $\Gamma'$-optimal groups.
\end{abstract}



\footnote{E-mail addresses: xuanlma@mail.bnu.edu.cn (X. Ma), fgmn\_1998@163.com (M. Feng), \\
wangks@bnu.edu.cn (K. Wang).}
\section{Introduction}
Each graph $\Gamma$ considered in this paper is a finite, simple and undirected graph with the vertex set $V(\Gamma)$ and the edge set $E(\Gamma)$.
Let $\Gamma_1$ and $\Gamma_2$ be graphs.  We call that $\Gamma_1$ is a {\em spanning subgraph} of $\Gamma_2$ if $V(\Gamma_1)=V(\Gamma_2)$ and $E(\Gamma_1)\subseteq E(\Gamma_2)$.
An {\em embedding}  from $\Gamma_1$ to $\Gamma_2$ is an injection $f:V(\Gamma_1)\rightarrow V(\Gamma_2)$ such that $\{x,y\}\in E(\Gamma_1)$ implies
$\{f(x),f(y)\}\in E(\Gamma_2)$. A graph $\Gamma_1$ can {\em embed} into a graph $\Gamma_2$ if there is an embedding from $\Gamma_1$ to $\Gamma_2$.

Let $G$ be a group.
The {\em power graph} $\Gamma_G$ of $G$ has the vertex set $G$ and two distinct elements are adjacent if one is a power of the other. In 2000,
Kelarev and  Quinn \cite{n1} introduced
the concept of a (directed) power graph. With this motivation, Chakrabarty, Ghosh and Sen \cite{CGS} introduced the undirected power graph of a group.
For convenience throughout we use the term ``power graph'' to refer to an undirected power graph defined as above. Recently, many interesting results on power graphs have been obtained, see \cite{AKC,AAA,Bu,acam,Cam,CGh,FMW,FMW1,kel21,kel2,kel22,MF,MFW}. A detailed list of results and open questions on power graphs can be found in \cite{AKC}.

Note that $\Gamma_G$ is complete if and only if $G$ is isomorphic to the cyclic group of prime power order (see \cite[Theorem 2.12]{CGS}).
This means that each graph can embed into $\Gamma_G$ for some finite group $G$. The motivation for this research is to decide the minimum order of a group $G$ such that a given graph can embed into $\Gamma_G$.


\begin{df}
The {\em power index} of a graph $\Gamma$,
denoted by $\Theta(\Gamma)$, is the least  order of a group $G$ such that  $\Gamma$ can embed into $\Gamma_G$.
Furthermore, this group $G$ is  {\em $\Gamma$-optimal} if $G$ has order $\Theta(\Gamma)$.
\end{df}

Let $K_n$ and $\mathbb Z_n$ be the complete graph and  cyclic group of order $n$, respectively. Since the power graph of any group of order $6$ is not complete, we have $\Theta(K_6)=\Theta(K_7)=7$, and $\mathbb{Z}_7$ is $K_6$-optimal and $K_7$-optimal.
For any graph $\Gamma$, one has
\begin{equation}\label{e11}
|V(\Gamma)|\le \Theta(\Gamma).
\end{equation}
Another motivation of this paper is to study which graphs  satisfy  the   equality in (\ref{e11}).

\begin{df}
A graph $\Gamma$ is {\em power-critical} if its order equals to $\Theta(\Gamma)$.
\end{df}

All power graphs are power-critical, but not vice versa. Actually,
given a graph $\Gamma$, determining whether $\Gamma$ is  power-critical is equivalent to determining whether $\Gamma$ is isomorphic to a spanning subgraph of the power graph of some finite group.
Many researchers \cite{CGh,FMW,mrs} investigated two groups which have isomorphic power graphs. Actually, this problem is to find two finite groups $G$ and $H$ such that $G$ is $\Gamma_H$-optimal and $H$ is $\Gamma_G$-optimal.

A cycle and a path of length $n$ are denoted by $C_n$ and $P_n$, respectively.
 It follows from \cite[Theorem~4.13]{CGS} that $\Gamma_{\mathbb{Z}_n}$ has a  Hamiltonian cycle, i.e., a cycle visiting each vertex of the graph . Thus, all paths and cycles are power-critical. Furthermore,
a group $G$ of order $n$ is $P_n$-optimal (resp. $C_n$-optimal) if and only if $\Gamma_G$ has a Hamilton path (resp. Hamilton cycle).

Let $K_{s,t}$ be the complete bipartite graph and $nK_2$ be the $1$-factor, i.e., the graph union of $n$ copies of $K_2$. This paper  focuses on the power indices of $K_n$, $K_{s,t}$ and $nK_2$. In Section 2, we compute $\Theta(K_n)$ and show that a $K_n$-optimal group is cyclic. In Section 3,  we give a necessary and sufficient condition for power-critical $K_{s,t}$, and under this condition,  all $K_{s,t}$-optimal groups are classified.
In Section 4, we show that $nK_2$ is power-critical, and give a necessary and sufficient condition for $nK_2$-optimal groups.

\section{Complete graphs}

Denote by $\rho_n$ the smallest prime power at least $n$.
Since $\rho_n$ is a prime power, any graph $\Gamma$ of order $n$ can  embed into  $\Gamma_{\mathbb{Z}_{\rho_n}}$. Hence, we have
$$
n=\Theta(\overline{K_n})\le\Theta(\Gamma)\le\Theta(K_n)\le\rho_n,
$$
where $\overline{K_n}$ is the null graph of order $n$. Moreover, any group of order $n$ is $\overline{K_n}$-optimal.

\begin{ob}\label{kn=n}
Let $n$ be a positive integer.

{\rm(i)} The complete graph $K_n$ is power-critical if and only if $n$ is a prime power.

{\rm(ii)} If $n$ is a prime power, then any graph of order $n$ is power-critical.
\end{ob}

In the following, we compute $\Theta(K_n)$ and find all groups which are $K_n$-optimal.
Let $n=p_1^{r_1}p_2^{r_2}\cdots p_m^{r_m}$,
where $p_1<p_2<\cdots<p_m$ are primes and $r_i\ge 1$
for $1\le i \le m$. Write
\begin{eqnarray*}
\chi_n&=&\varphi(n)+\varphi\left(\frac{n}{p_1}\right)+\cdots
+\varphi\left(\frac{n}{p_1^{r_1}}\right)+\varphi\left(\frac{n}{p_1^{r_1}p_2}\right)
+\varphi\left(\frac{n}{p_1^{r_1}p_2^2}\right)\\
& &+\cdots+\varphi\left(\frac{n}{p_1^{r_1}p_2^{r_2}}\right)+
\varphi\left(\frac{n}{p_1^{r_1}p_2^{r_2}p_3}\right)\\
& &+\cdots+\varphi\left(\frac{n}{p_1^{r_1}p_2^{r_2}\cdots p_m^{r_m-1}}\right)+\varphi(1),
\end{eqnarray*}
where $\varphi$ is the Euler totient function.
Recently, this number $\chi_n$ is studied by Curtina and Pourgholi \cite{cp}.

\begin{lem}{\rm\cite{cp}}\label{3xp}
  {\rm(i)} $\chi_n=\varphi(n)+\chi_{\frac{n}{p}}$, where $p$ is the least prime factor of $n$.

  {\rm(ii)} $\chi_n\le n$, with equality if and only if $n$ is a prime power.

  {\rm(iii)} $\chi_n=n-1$ if and only if $n$ is  twice an odd prime.
\end{lem}

Given a graph $\Gamma$, a subset  of $V(\Gamma)$ is a {\em clique} if its induced subgraph  is complete. The {\em clique number} of $\Gamma$, denoted by $\omega(\Gamma)$, is the maximum size of a clique in $\Gamma$.

\begin{lem}{\rm\cite{man}}\label{3lem}
  {\rm (i)} $\omega(\Gamma_{\mathbb{Z}_n})=\chi_n$.

  {\rm(ii)} For a finite group $G$, each clique  in $\Gamma_G$ is a clique in the power graph of a cyclic subgroup of $G$.
\end{lem}

\begin{thm}\label{cgraph}
For a positive integer $n$, we have
$$
\Theta(K_n)=\min\{k: \chi_k\ge n\}.
$$
Moreover,  a group $G$ is $K_n$-optimal if and only if $G$ is a cyclic group of order $\Theta(K_n)$.
\end{thm}
\proof
Write
 $$
 t=\min\{k: \chi_k\ge n\}.
 $$
Then $\omega(\Gamma_{\mathbb{Z}_t})=\chi_t$
by Lemma~\ref{3lem}~(i).
Note that $\chi_t\ge n$. Thus, there exists an embedding from $K_n$ to
$\Gamma_{\mathbb{Z}_t}$. This implies that
$\Theta(K_n)\le t$.

Suppose that $G$ is a $K_n$-optimal group. Then $\omega(\Gamma_{G})\ge n$ and $|G|=\Theta(K_n)$. By Lemma~\ref{3lem}, there exists an element $x\in G$ such that
$$
\chi_{|x|}=\omega(\Gamma_{\langle x\rangle})=\omega(\Gamma_{G})\ge n.
$$
It follows that
$$
|G|\ge|x|\ge t\ge \Theta(K_n)=|G|,
$$
 which implies that $\Theta(K_n)=t$ and $G$ is cyclic.

 Now suppose that $G$ is a cyclic group of order $t$. Then $\omega(\Gamma_G)=\chi_{t}\ge n$, and so $G$ is $K_n$-optimal.

 We accomplish our proof.
\qed

\medskip

By Observation~\ref{kn=n}, we have $n+1\le \Theta(K_n)\le\rho_n$ for a positive integer $n$ which is not a prime power.
We determine all $n$ satisfying $\Theta(K_n)=n+1$.

\begin{cor}\label{knn1}
Let $n$ be a positive integer which is not a prime power.
Then $\Theta(K_n)=n+1$ if and only if $n+1$ is a prime power or twice an odd prime.
\end{cor}
\proof
Assume that $\Theta(K_n)=n+1$. Then $n\le \chi_{n+1}\le n+1$.
If $\chi_{n+1}=n+1$, then $n+1$ is a prime power Lemma~\ref{3xp}~(ii). If $\chi_{n+1}=n$, then $n+1$ is twice an odd prime by Lemma~\ref{3xp}~(iii).

Conversely, we have $\chi_{n+1}=n+1$ or $n$. Note that $\chi_n<n$.
It follows that $\Theta(K_n)=n+1$, as desired.
\qed

\medskip

It follows from Theorem~\ref{cgraph} and Corollary~\ref{knn1} that $\Theta(K_{14})=16=\rho_{14}$.
In view of $\chi_{36}=27$, we get $\Theta(K_{34})=37=\rho_{34}$.
In view of $\chi_{93}=91$, one has $\Theta(K_{91})=93<\rho_{91}$. It is interesting to determine all  $n$ such that $\Theta(K_n)=\rho_n$.

\section{Complete bipartite graphs}

We begin this section by observing that the star $K_{1,t}$ can embed into the power graph of each group of order at least $1+t$.

\begin{ob}\label{}
The star $K_{1,t}$ is power-critical. Particularly, any group of order $1+t$ is $K_{1,t}$-optimal.
\end{ob}

In the remaining of this section, we always assume that  the complete bipartite graph $K_{s,t}$ have the vertex set partition $\{U,W\}$, and $2\le s\le t$. For a subset $S\subseteq V(K_{s,t})$ and an embedding
\begin{equation}\label{4f}
f: V(K_{s,t})\rightarrow V(\Gamma_G),
\end{equation}
write $f(S)=\{f(v):v\in S\}$.

\begin{lem}\label{subgroup}
Let $u$ and $v$ be vertices of $K_{s,t}$. With reference to {\rm(\ref{4f})},
if $|f(u)|=|f(v)|=p$ for some prime $p$, then $\langle f(u)\rangle=\langle f(v)\rangle$. Moreover, if $|G|=s+t$ and $p$ is a prime divisor of $s+t$, then $G$ has a unique
subgroup of order $p$.
\end{lem}
\proof
If $u$ and $v$ are adjacent, then one of $f(u)$ and $f(v)$ is a power of the other, and so $\langle f(u)\rangle=\langle f(v)\rangle$. Now suppose that $u$ and $v$ are nonadjacent. Without loss of generality, let $\{u,v\}\subseteq U$. Note that $2\le s\le t$. Pick $w\in W$ such that $f(w)$ is not the identity of $G$. Then both $f(u)$ and $f(v)$ are adjacent to $f(w)$ in $\Gamma_G$. Note that $|f(u)|=|f(v)|=p$. It is easy to check that $\langle f(u)\rangle=\langle f(v)\rangle$.

If $|G|=s+t$, then $\{f(U),f(W)\}$ is a partition of $G$. Hence, the desired result follows.
\qed

\begin{lem}\label{npgroup}
As refer to {\rm(\ref{4f})}, if $G$ is not a $p$-group and $|G|=s+t$, then $G$ is cyclic. In particular, one of $f(U)$ and $f(W)$ is a subset of $A$, where $A$ is the set of generators and the identity of $G$.
\end{lem}
\proof
Let $X$ be the set of all elements of prime order of $G$.
Note that $X$ contains two elements with distinct orders.
Then $X\subseteq f(U)$ or $f(W)$. Without loss of generality, let $X\subseteq f(U)$. Take any non-identity element $x$ in $f(W)$. Then each element of $X$ belongs to $\langle x\rangle$.

We now claim that every element of $f(U)$ belongs to $\langle x\rangle$.
In fact, suppose that there exists $y$ in $f(U)$ such that $y\notin \langle x\rangle$. Then $x\in\langle y\rangle$. Hence,
 there exist prime $p$ and positive integer $m$
such that $|y|$ is divisible by $p^m$ and $|x|$ is not
divisible by $p^m$. Let $z$ be an element of $\langle y\rangle$ of order $p^m$. Since $|x|$ is not a prime power,
$x$ and $z$ are not adjacent in $\Gamma_G$, which means that $z\in f(W)$.
Moreover, it is clear that there
exists an element $z'$ of order $q$ in $G$, where $q$ is a
prime divisor of $|G|$ and $q\ne p$.
Since $z'\in X\subseteq f(U)$, $z'$ is adjacent to $z$ in $\Gamma_G$, and so
$|z'|$ divides $|z|$ or $|z|$ divides $|z'|$, which contradicts
the orders of $z$ and $z'$. Thus, our claim is valid.

Take any non-identity element $x'$ in $f(W)\setminus\{x\}$, and let
$|x'|=p_1^{r_1}p_2^{r_2}\cdots p_l^{r_l}$, where $p_1,\ldots,p_l$ are pairwise distinct  primes and $r_i\ge 1$ for each $1\le i \le l$.
Let $u_i$ be an element of $\langle x'\rangle$ of order $p_i^{r_i}$.
Since there exists an element of $X$ of prime order which is different from $p_i$, one has $u_i\in f(U)$. By the claim above, we have $u_i\in \langle x\rangle$ for each $1\le i \le l$.
Since $\langle u_1\rangle\langle u_2\rangle\cdots\langle u_l\rangle=\langle x'\rangle$, one gets $x'\in \langle x\rangle$.
It follows that any element of $G$ belongs to $\langle x\rangle$, and
so $G=\langle x\rangle$, as wanted.
\qed

\begin{lem}\label{4optimal}
Let $2\le s\le t$.  If $\varphi(s+t)\ge s-1$, then $\mathbb{Z}_{s+t}$ is $K_{s,t}$-optimal.
\end{lem}
\proof Without loss of generality, assume that  $|U|=s$ and $|W|=t$.
Let $A$ be the set of generators and the identity of $\mathbb Z_{s+t}$.
Then $|A|=\varphi(s+t)+1$. Hence, there exists a bijection $h$ from $V(K_{s,t})$ to $\mathbb Z_{s+t}$ such that $h(U)\subseteq A$. It is easy to verify  that $h$ is an embedding from $K_{s,t}$ to $\Gamma_{\mathbb Z_{s+t}}$. Hence, the required result follows.
\qed

\begin{thm}\label{kmnl}
Let $2\le s\le t$.
 Then $K_{s,t}$ is power-critical if and only if
$\varphi(s+t)\ge s-1.$
\end{thm}
\proof
Suppose that $K_{s,t}$ is power-critical. With reference to (\ref{4f}), assume that $|G|=s+t$.
If $G$ is a $p$-group, then $s+t$  is a prime power,
which implies that
$$
\varphi(s+t)+1\ge \frac{s+t}{2}\ge s.
$$
If $G$ is not a $p$-group, it follows from Lemma \ref{npgroup} that $G$ is cyclic and
$$
s=\min\{|f(U)|,|f(W)|\}\le \varphi(s+t)+1.
$$

For the converse, the required result follows from Lemma~\ref{4optimal}.
\qed

\medskip

By Theorem~\ref{kmnl} and
some properties of the Euler's totient function, we get

\begin{cor}
{\rm(i)} For any positive integer $s$, there exists a positive integer $t_s$ such that $K_{s,t}$ is power-critical for each integer $t\ge t_s$.

{\rm(ii)} Let $s\ge 2$. Then
$K_{s,s}$ is power-critical if and only if $s$ is an odd prime or a power of $2$.
\end{cor}

For any $2\le s \le 6$ and $t\ge 7$, $K_{s,t}$ is power-critical. However, $\Theta(K_{6,6})=13$ and $\Theta(K_{9,9})=19$.

For $n\ge 2$, the generalized quaternion group $Q_{4n}$ is defined by
\begin{equation*}\label{2}
Q_{4n}=\langle x,y: x^{n}=y^2, x^{2n}=1, y^{-1}xy=x^{-1}\rangle.
\end{equation*}

\begin{cor}\label{op1}
Let $3\le k$ and $2\le s\le t$. Suppose that $K_{s,t}$ is power-critical.

$\rm{(i)}$ A group $G$ is $K_{2,2^k-2}$-optimal if and only if $G$ is isomorphic to $\mathbb{Z}_{2^k}$ or $Q_{2^k}$.

$\rm{(ii)}$ If $(s,t)\ne(2,2^k-2)$, then
a group $G$ is $K_{s,t}$-optimal if and only if
$G$ is isomorphic to $\mathbb{Z}_{s+t}$.
\end{cor}
\proof
(i) It is clear that $\mathbb{Z}_{2^k}$ is $K_{2,2^k-2}$-optimal.
Since $Q_{2^k}$ has a unique involution which is a power of any other nonidentity elements of $Q_{2^k}$, it is easy to check that $Q_{2^k}$ is $K_{2,2^k-2}$-optimal.
Now suppose that $G$ is $K_{2,2^k-2}$-optimal. Then $|G|=2^k$.
It follows from Lemma \ref{subgroup} that $G$ has a unique subgroup of order $2$, and so $G$ is isomorphic to $\mathbb{Z}_{2^k}$ or $Q_{2^k}$ by \cite[Theorem 5.4.10 (ii)]{Gor}.

(ii) It follows from Lemma~\ref{4optimal} that $\mathbb{Z}_{s+t}$ is $K_{s,t}$-optimal. In the following, suppose that $G$ is $K_{s,t}$-optimal. Then $|G|=s+t$.
If $G$ is not a $p$-group, then the required result holds by Lemma~\ref{npgroup}.
Now let $G$ be a $p$-group. Then $G$ is cyclic or generalized quaternion  by Lemma \ref{subgroup} and \cite[Theorem 5.4.10 (ii)]{Gor}.
Suppose that $G$ is isomorphic to $Q_{2^k}$ for some integer $k\ge 3$.  Since $(s,t)\ne (2,2^k-2)$ and $s+t=2^k$, we have $3\le s\le 2^{k-1}$. Hence, the power graph $\Gamma_G$ has $3$ distinct vertices with degree at least $2^{k-1}$, which contradicts the structure of $\Gamma_{Q_{2^k}}$ as shown in Figure~\ref{gq4n}.
\begin{figure}[hptb]
 \centering
  \includegraphics[width=8cm]{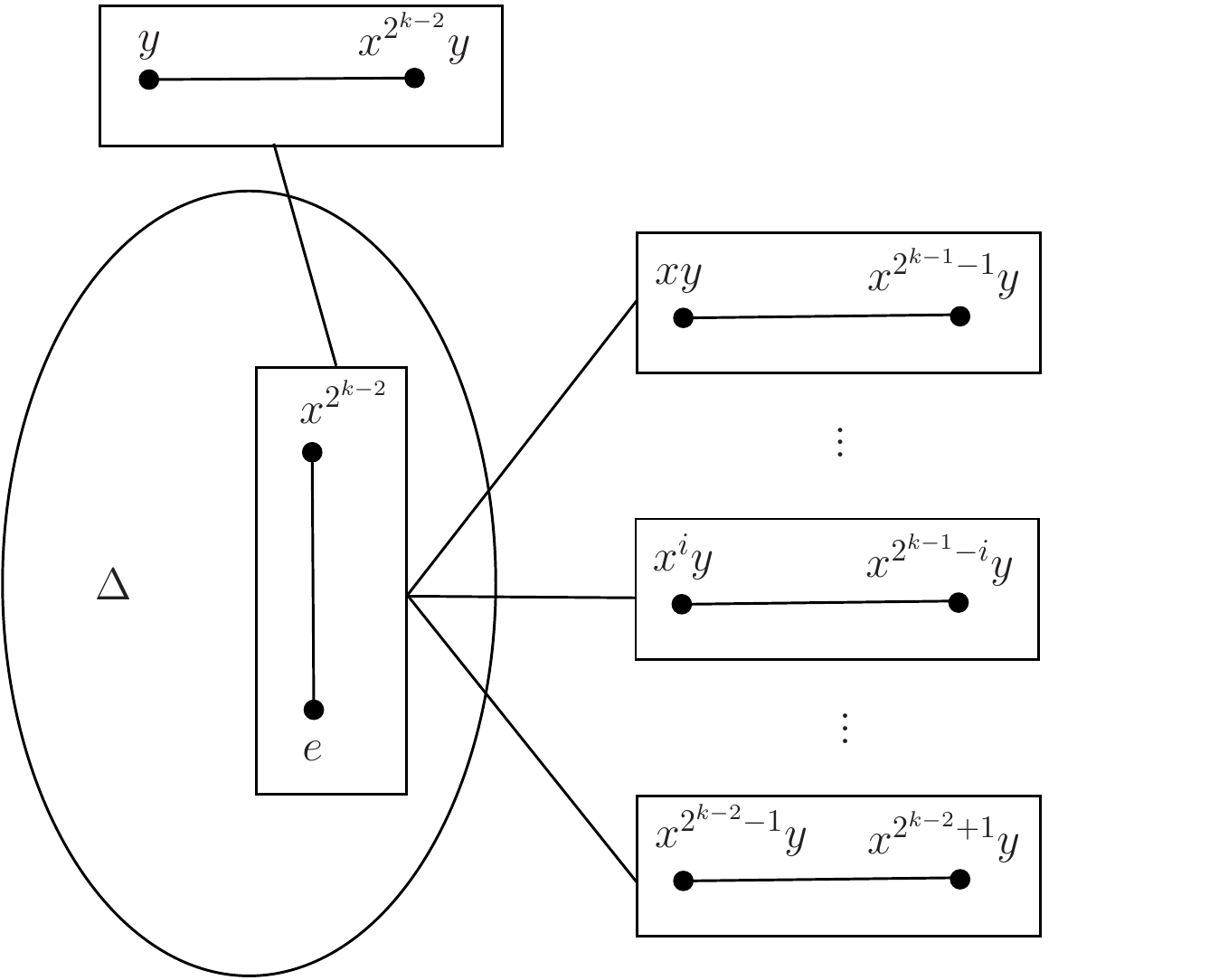}\\
  \caption{The power graph $\Gamma_{Q_{2^k}}$, where $e$ is the identity of $Q_{2^k}$ and  $\Delta$ is a complete graph of order $2^{k-2}$.}\label{gq4n}
\end{figure}

Therefore, the desired result follows.
\qed

\medskip

The following result is immediate by  Corollary~\ref{op1}.
\begin{cor}
Let $G$ be a group of order $n$. Then there exists a nonidentity element with degree at least $n-1$ in $\Gamma_G$  if and only if
$G$ is a cyclic group or  a generalized quaternion $2$-group.
\end{cor}

\section{$1$-factors}
Let $K_1+nK_2$ be the graph obtained from $nK_2$ by adding a vertex and joining this vertex to each vertex of $nK_2$.
In this section, we always use $e$ to denote the identity of a group $G$.

\begin{pro}\label{5k1nk2}
The graph $K_1+nK_2$ is power-critical, and any group $G$ of order $2n+1$ is $(K_1+nK_2)$-optimal.
\end{pro}
\proof We only need to show that any group $G$ of order $2n+1$ is $(K_1+nK_2)$-optimal. For each $x\in G\setminus\{e\}$,
we get $x\ne x^{-1}$. Write
\begin{eqnarray*}
 V(K_1+nK_2)&=&\{u,v_1,\ldots,v_n,w_1,\ldots,w_n\},\\
 E(K_1+nK_2)&=&\bigcup_{i=1}^n\{\{u,v_i\},\{u,w_i\},\{v_i,w_i\}\}.
\end{eqnarray*}
Then there exists a bijection $f$ from $V(K_1+nK_2)$ to $G$ such that
$f(u)=e$ and $f(v_i)=f(w_i)^{-1}$ for each $i\in\{1,\ldots,n\}$.
It is easy to verify that $f$ is an embedding  from $K_1+nK_2$ to $\Gamma_G$. Hence, the desired result follows.
\qed

\medskip

A {\em matching} in a graph  is a set of pairwise nonadjacent edges.
A {\em perfect matching} of a graph $\Gamma$ is a matching if its size equals to $\frac{|V(\Gamma)|}{2}$, and a {\em near-perfect matching} of  $\Gamma$ is a matching if its size equals to $\frac{|V(\Gamma)|-1}{2}$.
The following result is immediate by Proposition~\ref{5k1nk2}.

\begin{cor}
 If $G$ is a group of odd order, then $\Gamma_G$ has a near-perfect matching.
\end{cor}

 A path  $P$ in a power  graph $\Gamma_G$ is {\em inverse-closed }if $V(P)=\{x^{-1}:x\in V(P)\}$. Actually, for any distinct elements $u$ and $v$ in $G$, there exists an inverse-closed path  between $u$ and $v$ in $\Gamma_G$.
 Denote by $L(P)$ the set of the endpoints of a path $P$.

\begin{algorithm}[h!]
\caption{Find $x_1,\ldots,x_m$} \label{Palg}
\begin{algorithmic}[1]
\STATE {\bf Input}\qquad $u_1,\ldots,u_r$
\STATE Put $x_1=u_1$, $i=1$ and $m=1$
\FOR{$x_m\ne u_r$ and $x_m\ne u_r^{-1}$}
\STATE  Let $l_i$ be the index with $u_{l_i}=u_i^{-1}$. \label{p1}
\STATE $i\leftarrow\max\{i,l_i\}+1$ and $m\leftarrow m+1$ \label{p2}
\STATE $x_m=u_i$\label{p3}
\ENDFOR
\IF{$x_m=u_r$}
\STATE $x_m\leftarrow u_r^{-1}$
\ENDIF
\STATE {\bf Output}\qquad $x_1,\ldots,x_m$
\end{algorithmic}
\end{algorithm}

\begin{lem}\label{5lem}
   Let $P$ be an inverse-closed path in  a power graph $\Gamma_G$ such that  $|x|\ge 3$ for each $x\in V(P)$. Then there exists an inverse-closed path
  $$
  P'=(x_1,x_1^{-1},x_2,x_2^{-1},\ldots,x_m,x_m^{-1})
  $$
  such that $V(P')\subseteq V(P)$ and $L(P')=L(P)$.
\end{lem}
 \proof Write $P=(u_1,\ldots,u_r)$. By Algorithm~\ref{Palg}, we obtain some vertices $x_1,\ldots,x_m$ from $V(P)$. From Step~\ref{p1}, each vertex in $\{u_i,u_{l_i}\}$ is adjacent to $u_{i+1}$ and $u_{l_i+1}$ in $\Gamma_G$. It follows from Steps~\ref{p2} and~\ref{p3} that $\{x_j^{-1},x_{j+1}\}\in E(\Gamma_G)$ for $1\le j\le m-1$. Note that $\{x_1,\ldots,x_m\}\cap\{x_1^{-1},\ldots,x_m^{-1}\}=\emptyset$.
 Then $(x_1,x_1^{-1},x_2,x_2^{-1},\ldots, x_m,x_m^{-1})$ is the desired path.
 \qed

 \medskip

For a group $G$ of even order, we always use $U$ to denote the set of all involutions, and write $\overline U=U\cup\{e\}$.
Note that a group of even order has odd number of involutions.

\begin{thm}\label{perm}
  Let $G$ be a group of order $2n$ with $2k-1$ involutions. Then the followings are equivalent.

  {\rm(i)} The group $G$ is $nK_2$-optimal.

  {\rm(ii)} The power graph $\Gamma_G$ has a perfect matching.

  {\rm(iii)} There exist $k$ vertex-disjoint and inverse-closed paths $P_1,\ldots,P_{k}$ in $\Gamma_G$ such that
  $$
  \bigcup_{i=1}^k L(P_i)=\overline U.
  $$
\end{thm}
\proof It is clear that (i) and (ii) are equivalent. In the following, we shall show that (ii) and (iii) are equivalent.

Suppose (ii) holds. Let $\mathcal M$ be a perfect matching of $\Gamma_G$.
Using Algorithm~\ref{alg}, we obtain $k$ vertex-disjoint paths $P_1,\ldots,P_{k}$ such that $\bigcup_{i=1}^k L(P_i)=\overline U$. Hence (iii) holds.

Suppose (iii) holds. For $i\in\{1,\ldots,k\}$, write
$$
P_i=(u_{i1},x_{i1},x_{i2},\ldots,x_{it_i},u_{i2}).
$$
Without loss of generality, assume that $u_{k1}=e$. For $1\le i\le k-1$, we have $t_i\ge 2$. By Lemma~\ref{5lem},  there exists
$$P_i'=(u_{i1},y_{i1},y_{i1}^{-1},\ldots,y_{im_{i}},y_{im_i}^{-1},u_{i2})$$
such that $V(P_i')\subseteq V(P_i)$. Let $P_k'=(u_{k1},u_{k2})$.
Now we may assume that
$$
G\setminus\bigcup_{i=1}^kV(P_i')
=\{v_1,v_2,\ldots,v_r\}\cup \{v_1^{-1},v_2^{-1},\ldots,v_r^{-1}\},
$$
where $|v_i|\ge 3$ for $1\le i \le r$.
Then
$$
\begin{array}{l}
\big\{\{u_{11},y_{11}\},\{y_{11}^{-1},y_{12}\},\{y_{12}^{-1},y_{13}\},\ldots,
\{y_{1m_1}^{-1},u_{12}\},\\
\cdots \quad\cdots\quad\cdots\quad\cdots\quad
\cdots\quad\cdots\quad\cdots\quad\cdots\quad\cdots\\
\{u_{(k-1)1},y_{(k-1)1}\},\{y_{(k-1)1}^{-1},y_{(k-1)2}\},\{y_{(k-1)2}^{-1},
y_{(k-1)3}\},
\ldots,\{y_{(k-1)m_{k-1}}^{-1},u_{(k-1)2}\},\\
\{u_{k1},u_{k2}\},\{v_1,v_1^{-1}\},\ldots,\{v_r,v_r^{-1}\}\big\}
\end{array}
$$
is a perfect matching  of $\Gamma_G$.
\qed

\begin{algorithm}[h!]
\caption{Find $k$ vertex-disjoint paths $P_1,\ldots,P_{k}$ such that $\bigcup_{i=1}^k L(P_i)=\overline U$} \label{alg}
\begin{algorithmic}[1]
\STATE {\bf Input}\qquad $\overline U$ and $\mathcal M$
\STATE  Set $A:=\overline U$\qquad // $|A|=2k$
\FOR{$i=1,\ldots,k$}
\STATE Choose $u\in A$
\STATE Put $P_i:=(u,x)$ for $\{u,x\}\in\mathcal M$
\WHILE{$x\notin A$}
\STATE $P_i\leftarrow(P_i,x^{-1},y)$ for $\{x^{-1},y\}\in\mathcal M$\qquad // Note that $x\notin\overline U\Leftrightarrow x\ne x^{-1}$
\STATE $x\leftarrow y$
\ENDWHILE
\STATE $A\leftarrow A\setminus\{u,x\}$
\ENDFOR \qquad//$A=\emptyset$
\STATE {\bf Output}\qquad $P_1,\ldots,P_k$
\end{algorithmic}
\end{algorithm}

Note that $\mathbb{Z}_{2n}$ has a unique involution, which is adjacent to the identity in $\Gamma_{\mathbb{Z}_{2n}}$. The following result holds from Theorem~\ref{perm}.

\begin{thm}\label{nk2}
The $1$-factor is power-critical and $\mathbb{Z}_{2n}$ is $nK_2$-optimal.
\end{thm}

Actually, any group of order $2n$ having a unique involution is
$nK_2$-optimal. Hence, $Q_{4n}$ is $2nK_2$-optimal.
Let $D_{2n}$ be the dihedral group of order $2n$. In $\Gamma_{D_{2n}}$, any path between two distinct involutions contains the identity, so $D_{2n}$ is not $nK_2$-optimal.

\section*{Acknowledgement}
This work was carried out during Ma Xuanlong's visit to the Beijing Normal University(Dec. 2016--Jan. 2017).
K. Wang's research was supported by National Natural Science Foundation of China (11371204, 11671043) and the Fundamental Research Funds for the Central University of China.

\end{document}